\documentclass[12pt, oneside,reqno]{amsart}   	
\usepackage{geometry}
\usepackage{float}
\allowdisplaybreaks
\usepackage{graphicx}				
\usepackage{amssymb}
\usepackage{amsmath}
\usepackage[all]{xy}
\usepackage{mathtools}
\usepackage{tikz-cd}
\usepackage{enumitem}
\usepackage{tikz-cd}
\usepackage{appendix}

\usepackage[%
bookmarks=true,
colorlinks,
linkcolor=blue,
urlcolor=blue,
citecolor=blue,
plainpages=false,
pdfpagelabels,
final,
breaklinks=true\between
]{hyperref}
\usepackage{hyperref}
\usepackage[numbers,sort&compress]{natbib}

\newtheorem{thm}{Theorem}[section]
\newtheorem{lem}[thm]{Lemma}

\newtheorem{conj}{Conjecture}[section]

\newtheorem{rem}{Remark}[section]

\numberwithin{equation}{section}

\def\Pb{\ifmmode{\Bbb P}\else{$\Bbb P$}\fi}
\def\Z{\ifmmode{\Bbb Z}\else{$\Bbb Z$}\fi}
\def\N{\ifmmode{\Bbb N}\else{$\Bbb N$}\fi}
\def\C{\ifmmode{\Bbb C}\else{$\Bbb C$}\fi}
\def\R{\ifmmode{\Bbb R}\else{$\Bbb R$}\fi}
\def\S{\ifmmode{S^2}\else{$S^2$}\fi}

\def\S{\cal S}

\newenvironment{pf}{\paragraph{Proof:}}{\hfill$\square$ \newline}

\begin{document}

\title[Ancient doughnuts]{On the construction of closed nonconvex nonsolition ancient mean curvature flows}
\author{Theodora Bourni, Mathew Langford,  and Alexander Mramor   }
\address{Department of Mathematics, University of Tennessee Knoxville,
Knoxville TN, 37996}
\email{ tbourni@utk.edu,  mlangford@utk.edu}
\address{Department of Mathematics, Johns Hopkins University, Baltimore, MD, 21231}
\email{amramor1@jhu.edu}

\begin{abstract} We construct closed, embedded, ancient mean curvature flows in each dimension $n\ge 2$ with the topology of $S^1\times S^{n-1}$. These examples are not mean convex and not solitons. They are constructed by analyzing perturbations of the self-shrinking doughnuts constructed by Drugan and Nguyen (or, alternatively, Angenent's self shrinking torus when $n =2$)
\end{abstract}
\maketitle
\section{Introduction}

The mean curvature flow, whereby a submanifold is deformed with velocity equal to its mean curvature vector, is a natural analogue of the heat equation in submanifold geometry, and has a multitude of applications in topology, geometry, physics and materials science. \emph{Ancient} solutions (those solutions defined for times $t \in (-\infty, T)$, $T \leq \infty$) are of particular interest due to their natural role in modeling singularity formation (via `blow-up' analysis). They are also interesting since they are the parabolic counterpart of proper minimal surfaces; indeed given the general ``improving'' effect of heat equations they should have special properties which distinguish them from a typical flow. 
$\medskip$

A special class of ancient solutions are the solitons. These solutions evolve by a combination of rigid motions and dilations, which reduces the mean curvature flow equation to an elliptic equation for a given timeslice. In simple cases, the ancient solutions which arise as blow-ups of singularities can be shown to be soliton solutions, but this is not necessarily the case in general. In this paper we construct new such ancient flows exhibiting new phenomena, as we elaborate below: 
$\medskip$

\begin{thm}\label{thm:ancient doughnut}
For each $n\ge 2$, there exists a closed, embedded, nonsoliton ancient mean curvature flow $\{M_t\}_{t\in(-\infty,0)}$ in $\R^{n+1}$ with the topology of $S^1\times S^{n-1}$. 
\end{thm}
$\medskip$

It will be evident from our construction that $\{M_t\}_{t\in(-\infty,0)}$ is rotationally symmetric and becomes singular in a circle of neckpinch singularities at time $0$. So it must become mean convex before the final time. On the other hand, since its blowdown, $\lim_{\lambda\to 0}\{\lambda M_{\lambda^{-2}t}\}_{t\in(-\infty,0)}$, is a self-shrinking $S^1\times S^{n-1}$, $M_t$ cannot be mean convex for all times. However, our solution does satisfy a natural mean convexity condition related to self-shrinker mean convexity (namely, inequality \eqref{eq:shrinker mean convex parabolic} below\footnote{Mean curvature flows satisfying this condition were studied by Smoczyk \cite{Smk} and Lin \cite{Lin}.}).
$\medskip$

The simplest nonsoliton example is the Angenent oval solution to the curve shortening flow $\{\Gamma_t\}_{t\in(-\infty,0)}$, where $\Gamma_t:=\{(x,y)\in \R^2:\cos x =\mathrm{e}^t\cosh y\}$ \cite{Ang92}. In higher dimensions, there is a family of compact ``ancient ovaloids'' interpolating between a cylinder and a sphere \cite{HaHe,Wa11,Wh00} and also compact ``ancient pancakes'' which lie in slab regions \cite{BLT1,Wa11}. These examples are all convex. Indeed, under certain `noncollapsing' conditions, ancient solutions are necessarily convex \cite{HK1,L17,Lin}.
$\medskip$

On the other hand, the third author, with Payne, has constructed non-convex ancient solutions \cite{MramorPayne}. These examples are embedded and mean convex (but not convex), and evolve out of minimal surfaces at time $t=-\infty$. Since they evolve out of minimal surfaces in Euclidean space, they are necessarily noncompact --- the examples constructed presently though are compact embedded such flows. 
$\medskip$

For the curve shortening flow, Angenent and You \cite{AngenentYou} have constructed a very large family of non-convex solutions by gluing together an essentially arbitrary family of Grim Reapers along their common asymptotes. This family includes compact immersed examples and noncompact embedded examples, but no compact, embedded examples. Conjecturally there is also the ancient ``yin-yang spiral'' curve shortening flow, but, to our knowledge, this has yet to be rigorously constructed. 
$\medskip$

These examples show that nonconvex ancient solutions can be quite flexible. On the other hand, \emph{convex} ancient solutions appear to be quite rigid: the only convex ancient solutions to curve shortening flow are the stationary lines, shrinking circles, Grim Reapers and Angenent ovals \cite{DHS,BLT3} and in higher dimensions there are further partial results in this direction \cite{ADS1,ADS2,BLT1,BrendleChoi1,BrendleChoi2,HuSi15,KevinLiouville}. Moreover, an important result of X.-J.~Wang states that a convex ancient solution which does not sweep out all of $\R^{n+1}$ must sweep out a slab region \cite{Wa11}. In particular, no convex ancient solution sweeps out a half-space (cf. \cite{ChiniMollerAncient,ChiniMollerTranslators}).
$\medskip$

\section{Preliminaries} 

In this section we collect some general facts and observations about mean curvature flow which we will employ below. 
$\medskip$

An important consequence of the maximum principle is the comparison principle (also known as the avoidance principle), which says that two initially disjoint, proper hypersurfaces will remain disjoint over their common interval of existence if at least one of them is compact. In particular, by comparison with a sufficiently large enveloping sphere, every compact hypersurface must develop a singularity in finite time. 
$\medskip$

Huisken's monotonicity formula \cite{Hu90} implies that singularities are modeled, in a precise sense, by \emph{self-similarly shrinking solutions}. That is, solutions $\{\Sigma_t\}_{t\in(-\infty,0)}$ of the form $\Sigma_t=\sqrt{-t}\Sigma_{-1}$. The time $t=-1$ slice $\Sigma_{-1}$ of a self-similarly shrinking solution necessarily satisfies the self-shrinker equation,
\begin{equation}\label{eq:self-shrinker}
H(x)-\frac{\langle x, \nu(x) \rangle}{2} = 0\,,
\end{equation} 
where $H(x)$ denotes the mean curvature at $x\in\Sigma_{-1}$ corresponding to the choice of unit normal $\nu(x)$. Important below, these are also minimal surfaces in the Gaussian metric $G = e^{-|x|^2/4} g_{\mathrm{flat}}$. The related \textit{Gaussian area} is then given by
\begin{equation}
A(M) := \int_M e^{-\frac{|x|^2}{4}} d\mathrm{vol}
\end{equation}

Solutions to \eqref{eq:self-shrinker} are called self-shrinkers. Conversely, self-shrinkers $\Sigma$ give rise to ancient solutions $\{\sqrt{-t}\Sigma\}_{t\in(-\infty,0)}$ to mean curvature flow. Sometimes one is led to choose different blowup procedures though, which will in general only yield ancient mean curvature flows. 
$\medskip$

An important example for our construction below is the self-shrinking $S^1\times S^{n-1}$ constructed by\footnote{(Rotationally symmetric) self-shrinking $S^1\times S^{n-1}$ were first constructed by Angenent using \textsc{ode} methods \cite{Ang92}. It is believed that there is only one rotationally symmetric example. Since this is not yet known \cite{DrLeeNgu}, we use the example of Drugan and Nguyen, which as we discuss below has entropy less than 2.} Drugan and Nguyen in \cite{DrNgu}. It is rotationally symmetric and has entropy less than 2 \cite{DrNgu,YBK} (by entropy here we mean in the sense of Colding and Minicozzi \cite{ColdingMinicozziGeneric}: the supremum of recentered and rescaled Gaussian densities) --- this is explained more below in the proof. We will exploit this bound along with the fact that entropy is non-increasing under the flow and lower semicontinuous. 
$\medskip$

Generally speaking, the only satisfactory understanding of self-shrinkers to the mean curvature flow is in the mean convex case, where they are known (under some weak conditions) to be precisely spheres and (generalized) cylinders \cite{Hu90,Hu93,ColdingMinicozziGeneric}. We will analyze perturbations of self-shrinking $S^1\times S^{n-1}$'s which are \textit{self-shrinker mean convex}, that is,
\begin{equation}\label{eq:shrinker mean convex}
H-\frac{\langle x, \nu \rangle}{2} > 0.
\end{equation} 
This condition is neither scale invariant nor preserved under the flow; however, for flows defined on negative time intervals, the related condition
\begin{equation}\label{eq:shrinker mean convex parabolic}
H-\frac{\langle x, \nu \rangle}{-2t} > 0
\end{equation}
is both invariant under parabolic scaling and preserved under mean curvature flow. 
To see that it is preserved, observe that the function $F:=\langle x,\nu\rangle+2tH$ satisfies the Jacobi equation,
\[
\partial_tF = \Delta F + |A|^2F\,.
\]
This is easily proved directly \cite{Smk}. It also follows from the fact that it is the normal component of the variation field generated by parabolic dilations.
$\medskip$

As Andrews observed in \cite{An12}, two-sided noncollapsing is preserved with respect to $F$, in the sense that, for a compact, embedded solution $\{M_t\}_{t\in[t_0,T)}$ to mean curvature flow,
\begin{equation}\label{eq:Fnoncollapsing}
\min_{M_{t_0}}\frac{\underline k}{F}\le \min_{M_t}\frac{\underline k}{F}\;\;\text{and}\;\;\max_{M_t}\frac{\overline k}{F}\le \max_{M_{t_0}}\frac{\underline k}{F}
\end{equation}
for all $t>t_0$, where, at each $x\in M_t$, $\overline k$ (resp. $\underline k$) is the inverse radius of the largest oriented\footnote{Note that $\underline k$ can be negative.} ball which lies inside (resp. outside) of $M_t$ and touches it at $x$. We will refer to such flows as $F\textit{-noncollapsing}$. We note that Lin \cite{Lin} has extended the Haslhofer--Kleiner local curvature estimates \cite{HK1} to this setting.

\section{Proof of Theorem \ref{thm:ancient doughnut}} 

Denote by $T = T^n$ the time $t=-1$ slice of the self-shrinking $S^1\times S^{n-1}$ of Drugan and Nguyen. We will construct a sequence of hypersurfaces $T_i$ and corresponding flows $\{T^i_t\}_{t \in [-1, t_i)}$, $t_i<0$, with the following properties:
\begin{enumerate}
\item $T_i$ is rotationally symmetric and self-shrinker mean convex, and $T_i\to T$ in the smooth topology as $i \to \infty$.
\item 
For each $i$, $t_i$ is the singular time of $\{T^i_t\}_{t \in [-1, t_i)}$, $\lim_{i\to\infty} t_i=0$, and the corresponding limit set $\mathcal{S}_i$ is a round circle (whose radius we denote by $d_i>0$).
\item $d_i^{-2}t_i\sim -1$. In particular, $d_i \to 0$ as $i \to \infty$.
\item There exists $C<\infty$ such that, after parabolically rescaling by $d_i$ to obtain flows $\{\widetilde{T}^i_t\}_{t\in[-d_i^{-2},\tilde t_i)}$, where $\tilde T_t^i:=d_i^{-1}T^i_{d_i^{2}t}$ and $\tilde t_i:=d_i^{-2}t_i$, 
\[
(\tilde t_i-t)|\tilde A^i_{(p,t)}|^2 < C
\]
for each $i\in\N$, $t\in[-1,\tilde t_i)$ and $p\in \tilde T{}^i_t$, where $\tilde A^i$ is the second fundamental form of $\{\widetilde{T}^i_t\}_{t\in[-d_i^{-2},\tilde t_i)}$.
\end{enumerate} 
By items (3) and (4), a subsequence of the flows $\{\widetilde{T}^i_t\}_{t\in[-d_i^2, \tilde t_i)}$ converges in the smooth topology (uniformly on compact time subintervals) to an ancient mean curvature flow $\{M_t\}_{t\in(-\infty,t_\ast)}$. By (2) and (3), the limit flow is nonempty and its limit set (at time $t_\ast$) is a circle of radius 1. Moreover, by (3), the diameter is type I, in the sense of Lemma \ref{lem:radius} , and thus the limit flow is compact. Therefore, by item (1), its blowdown is a self-shrinking torus. 
$\medskip$

No timeslice $M_t$ is convex, although they do become mean convex at some time $t_0<t_\ast$. On the other hand, the solution $\{M_t\}_{t\in(-\infty,t_\ast)}$ does satisfy the inequality \eqref{eq:shrinker mean convex parabolic}. In particular, $\{M_t\}_{t\in(-\infty,t_\ast)}$ is not self-similar (this also follows from the fact that the singular set is a non-trivial circle).
$\medskip$

First we describe the construction of the hypersurface $T_i$ (and prove item (1)), which is straightforward. 

\begin{lem} For any $i$, there is a perturbation $T_i$ of $T$ which is self-shrinker mean convex, rotationally symmetric, and satisfies $\Vert T - T_i\Vert_{C^2} < C/i$, where $C< \infty$ is some constant depending on $T$.  
\end{lem}
\begin{pf}

For $i$ sufficiently large, let $T_i$ be the constant inward variation $-\frac{1}{i}\nu$ of $T$, where $\nu$ is the outward pointing unit normal to $T$. Then $T_i$ is automatically rotationally symmetric. Now we need to check that such perturbation will be self-shrinker mean convex. 
$\medskip$

To achieve with this, we recall (see \cite{CIMW}, Lemma 1.2) that, varying a self-shrinker in the direction $u \nu$, we have 
\begin{equation}
\left.\frac{d}{ds}\right|_{s = 0} \left(H - \frac{\langle x,\nu\rangle}{2}\right) = - Lu\,,
\end{equation} 
where $L$ is the Jacobi operator, $\Delta - \langle \frac{x}{2}, \nabla \cdot \rangle + |A|^2 + \frac{1}{2}$. It follows that small constant inward variations are self-shrinker mean convex. The $C^2$ estimate follows from the bounded geometry of the torus and the choice of perturbation.
\end{pf}
$\medskip$

For each $i\in\N$, let $\{T^i_t\}_{t\in[-1,t_i)}$ be the maximal smooth solution to mean curvature flow with initial condition $T^i_0=T_i$. Since each initial datum $T_i$ is self-shrinker mean convex, each solution $\{T^i_t\}_{t\in[-1,t_i)}$ satisfies \eqref{eq:shrinker mean convex parabolic}. Since, by the avoidance principle, $T^i_t$ is enclosed by the time $t$ slice $T_t$ of the self-shrinking $S^1\times S^{n-1}$, we have $t_i<0$ for each $i$. On the other hand, by the continuous dependence of solutions on initial data, $\liminf_{i\to\infty}t_i\ge 0$ (cf. \cite{KevinMaximalTime}). Thus, $t_i\to 0$. The convergence to a circle will follow from Lemma \ref{lem:circle limit} below, completing the proof of item (2). The control on its radius (item (3)) will then follow from the following lemma.

\begin{lem}\label{lem:radius}
Denoting by $\overline{d}_i(t)$ and $\underline{d}_i(t)$ the maximal and minimal distance of a point $T^i_t$ to the origin respectively, there exists $0 < C < \infty$, independent of $i$, so that
\begin{equation}
 \frac{1}{C} \sqrt{-2t} < \underline{d}_i(t) < \overline{d}_i(t) < C\sqrt{-2t}\,.
 \end{equation}
\end{lem} 
\begin{pf} 
This follows immediately from the fact that $T^i_t$ is enclosed by the time $t$ slice of the self-shrinking $S^1\times S^{n-1}$ (a consequence of the avoidance principle).
\end{pf} 
$\medskip$

Because of its importance below, we explain the estimate for the entropy of $T$ given by Drugan and Nguyen \cite{DrNgu}:
\begin{lem}\label{lem:entropy} The entropy of $T$ is strictly less than 2 for each $n$.  
\end{lem}
\begin{pf}
In \cite{DrNgu} (as in Angenent \cite{Ang92}), Drugan and Nguyen reduce the problem of finding rotationally symmetric self-shrinking $S^1\times S^{n-1}$'s to finding closed geodesics $\gamma$ for the metric
\begin{equation} 
g = \lambda^2 g_{\mathrm{flat}} = \lambda^2(dr^2 + dx^2), \text{ } \lambda = r^{n-1} e^{-\frac{1}{4}(x^2 + r^2)}
\end{equation}
defined in the half plane $\R^2_{+} = \{(r,x) \mid r > 0, x \in \R \}$. Denoting by $L_n(C)$ the length of a curve $C$ in this metric, the geodesics  found in \cite{DrNgu}, denoted $\gamma_\infty$, satisfy the following estimate \cite[Theorem 1]{DrNgu}:   
\begin{equation}
L_n(\gamma_\infty) < 2 \int\limits_0^\infty s^{n-1} e^{-s^2/4} ds = 2^n \Gamma\left(\frac{n}{2}\right).
\end{equation} 

As discussed in \cite{Ang92} (to be precise, see the equation after (5) in \cite{Ang92}), the corresponding Gaussian area $A(T)$ of the corresponding self shrinking $S^1\times S^{n-1}$ is then equal to $\mathrm{Vol}(S^{n-1}) L_n(\gamma) < 2^n n \pi^{\frac{n}{2}}  \frac{\Gamma(\frac{n}{2})}{\Gamma(\frac{n}{2} + 1)}$. Now we recall from \cite{ColdingMinicozziGeneric} the well known fact that the entropy of a compact self-shrinker is equal to the $F$ functional $F_{0,1}$, which is simply the Gaussian area normalized so that the plane has value 1. Thus,
\begin{equation} 
\lambda(T) = F_{0,1}(T) = \frac{1}{(4\pi)^{n/2}}A(T) < \frac{2^n n \pi^{\frac{n}{2}}}{(4\pi)^{n/2}} \frac{\Gamma(\frac{n}{2})}{\Gamma(\frac{n}{2} + 1)} = n \frac{\Gamma(\frac{n}{2})}{\Gamma(\frac{n}{2} + 1)}  = 2
\end{equation} 
Which is the bound we claimed. 
\end{pf}
\begin{rem}
A very good estimation of the entropy (1.85122) of the Angenent torus was recently obtained by Berchenko-Kogan \cite{YBK}, and its use when $n=2$ in lieu of Lemma \ref{lem:entropy} is valid throughout (in case they are indeed distinct).
\end{rem}

As the core of the argument, we prove the uniform curvature estimate of item (4). To achieve this, we exploit the rotational symmetry and low entropy of $T$.

\begin{lem}\label{lem:curvature bound}
There exists $C < \infty$, independent of $i$, so that 
\begin{equation} 
|A^i_{(p,t)}| < \frac{C}{\sqrt{t_i-t}}
\end{equation} 
for all $p\in M^i_t$ and $t \in [-1, t_i)$, where $A^i_{(\cdot,t)}$ is the second fundamental form of $M^i_t$.
\end{lem}  
\begin{pf}
Suppose there is in fact no such $C$. Then, following Hamilton (as in \cite[Section 4]{HuSi99a}), we can (after passing to a subsequence if necessary) choose, for each $i$, a time $s_i\in[-1,t_i)$ and a point $p_i\in T^i_{s_i}$ so that
\[
\left(t_i-\frac{1}{i}-s_i\right)\vert A^i_{(p_i,s_i)}\vert^2=\sup_{p\in T^i_{t},t\in [-1,t_i-\frac{1}{i}]}\left(t_i-\frac{1}{i}-t\right)\vert A^i_{(p,t)}\vert^2
\]
and
\[
\left(t_i-s_i\right)\vert A^i_{(p_i,s_i)}\vert^2\to\infty.
\]
Set $\lambda_i:=\vert A^i_{(p_i,s_i)}\vert$ and consider the flows $\{\tilde T{}^i_t\}_{t\in[-\lambda_i^{2}(1+s_i),\lambda_i^2(t_i-s_i))}$ defined by
\[
\tilde T{}^i_t:=\lambda_i\left(T^i_{s_i+\lambda_i^{-2}t}-p_i\right).
\]
After passing to a subsequence, the flows $\{\tilde T{}^i_t\}_{t\in[-\lambda_i^{2}(1+s_i),\lambda_i^2(t_i-s_i))}$ converge locally uniformly in the smooth topology to a non-flat eternal limit flow $\{\tilde T{}^\infty_t\}_{t\in(-\infty,\infty)}$ which, by the rotational symmetry of the sequence, the type II blowup rate, and Lemma \ref{lem:radius}, splits off an $(n-1)$-dimensional plane. We claim that the corresponding solution to curve shortening flow is convex. To see this, we make use of the self-shrinker mean convexity of the sequence to estimate, for any $\lambda_i(p-p_i)\in \tilde M^i_t=\lambda_i(M_{s_i+\lambda_i^{-2}t}-p_i)$,
\[
\tilde H_i(\lambda_i(p-p_i),t)=\lambda_i^{-1}H_i(p,s_i+\lambda_i^{-2}t)\ge \lambda_i^{-1}\frac{\langle p,\nu_i(p,s_i+\lambda_i^{-2}t)\rangle}{-2(s_i+\lambda_i^{-2}t)}\,.
\]
By Lemma \ref{lem:radius},
\[
\vert p\vert\le C\sqrt{-2(s_i+\lambda_i^{-2}t)}\,,
\]
where $C$ is a constant that depends only on $T$, and hence
\[
\tilde H_i(p,t)\ge \frac{-C}{\sqrt{-2(\lambda_i^2s_i+t)}}=\frac{-C}{\sqrt{2\big(\lambda_i^2(t_i-s_i)-\lambda_i^2t_i-t\big)}}\,.
\]
Since $t_i< 0$ and, by hypothesis, $\lambda_i^2(t_i-s_i)\to\infty$, we conclude that the limit flow $\{\tilde T{}^\infty_t\}_{t\in(-\infty,\infty)}$ is mean convex. It follows that the cross section is a non-flat, non-compact, convex ancient solution to curve shortening flow, which we conclude must be the Grim Reaper by the classification in \cite{BLT3}.
$\medskip$

Recall that entropy is invariant under translations and dilations, lower semicontinuous under taking limits, and monotone under the mean curvature flow. With that in mind, since the Grim hyperplane has entropy 2, whereas $T$ (and hence $T^i$ for $i$ large enough, since entropy is known to be continuous under compactly supported smooth perturbations) has entropy strictly less than 2, we arrive at a contradiction. 
\end{pf} 

The uniform curvature estimate in item (4) follows.

\begin{rem}
By compactness and self-shrinker mean convexity of the initial data $T_i$, and the maximum principle, each solution $\{T^i_t\}_{t\in[-1,t_i)}$ is $F$-noncollapsing \cite{An12}. This, and the rotational symmetry of the solutions, is sufficient to rule out type II singularities (using work of Lin \cite{Lin}) and conclude that the profile curve shrinks to a point, for a fixed $i$. However, since $T$ is a self-shrinker, the quality of noncollapsing degenerates as $i \to \infty$. So this argument does not provide the \emph{uniform} type I curvature estimate we seek in item (4) above.
\end{rem}

We now show that the profile curves shrink to a point, which proves the remaining claim of item (2).

\begin{lem}\label{lem:circle limit}
$T_t^i\to \mathcal{S}_i$ as $t\to t_i$, where $\mathcal{S}_i$ is a round circle.
\end{lem}
\begin{pf}

By the preceding lemma, $\{T^i_t\}_{t\in[-1,t_i)}$ is of type I. Since it is $F$-noncollapsing and rotationally symmetric, the convexity estimate of Lin\footnote{Note that Lin's local curvature estimate, and the resulting convexity estimate, only require the self-shrinker mean convexity condition \eqref{eq:shrinker mean convex parabolic}, rather than starshapedness. In dimension 2, we could also have applied the convexity estimate of Smoczyk \cite{Smk}.} \cite{Lin} implies that the tangent flow to $\{T^i_t\}_{t\in[-1,t_i)}$ about a point $p\in \mathcal{S}_i$ must be a shrinking cylinder (cf. \cite{KevinLiouville}). It follows that the profile curve of $\{T^i_t\}_{t\in[-1,t_i)}$ shrinks to a point, which implies the claim.
\end{pf} 
$\medskip$

Theorem \ref{thm:ancient doughnut} now follows by combining properties (1)-(4).

\section{Concluding remarks}

It seems reasonable to expect that the blowdown of the ancient solution above, in light of (1), is the self-shrinking $S^1\times S^{n-1}$ we started with. Of course this will certainly be true if the Angenent--Drugan--Nguyen self-shrinking $S^1\times S^{n-1}$ is unique amongst rotationally symmetric self-shrinking $S^1\times S^{n-1}$'s. 
$\medskip$

The construction explicitly required the rotational symmetry, topology, compactness, and low entropy of the perturbed self-shrinker. But we feel that our construction fits into a general phenomenon.

\begin{conj} Let $\Sigma$ be a compact, embedded, nonround self-shrinker in $\R^{n+1}$. Then there exists an associated nonsoliton ancient flow $\{M_t\}_{t\in(-\infty,0)}$. 
\end{conj} 
Part of the difficulty in establishing this conjecture following the approach above are the following two points we wish to emphasize: 
\begin{enumerate}
\item In general, there will be translators which have lower entropy than the initial self-shrinker --- and even possibly compact self-shrinkers with entropy greater than 2 (the entropy of the Grim Reaper).
\item Without the rotational symmetry hypothesis, our convexity argument in the proof of Lemma \ref{lem:curvature bound} only implies mean convexity of the blow-up limit (rather than convexity).
\end{enumerate} 
These are relevant because, in the curvature estimates established above, we used the entropy bound on $T$ along with the classification of convex ancient curve shortening flows to rule out type II curvature blowup.
$\medskip$
$\medskip$
$\medskip$
$\medskip$

\bibliographystyle{acm}
\bibliography{bibliography}

\end{document}